\documentclass[12pt,twoside]{amsart}
\usepackage{amssymb}
\usepackage{graphicx}

\nonstopmode

\numberwithin{equation}{section}
\font\tengothic=eufm10 scaled\magstep 1
\font\sevengothic=eufm7 scaled\magstep 1
\newfam\gothicfam
          \textfont\gothicfam=\tengothic
          \scriptfont\gothicfam=\sevengothic

\DeclareMathOperator{\pnt}{\raise 0.5mm \hbox{\large\bf.}}

\newtheorem{theorem}{Theorem}[section]
\newtheorem{lemma}[theorem]{Lemma}
\newtheorem{proposition}[theorem]{Proposition}
\newtheorem{corollary}[theorem]{Corollary}

\theoremstyle{definition}
\newtheorem{definition}[theorem]{Definition} 
\newtheorem{remark}[theorem]{Remark}
\newtheorem{example}[theorem]{Example}

\newtheorem{question}[theorem]{Question}

\begin{document}
\title{The strength of the Weak Lefschetz Property}

\author[Juan Migliore]{Juan Migliore${}^*$}
\address{
Department of Mathematics, University of Notre Dame, Notre Dame, IN
46556, USA}
\email{Juan.C.Migliore.1@nd.edu}

\author[Fabrizio Zanello]{Fabrizio Zanello}
\address{Department of Mathematical Sciences, Michigan Technological University, Houghton, MI 49931, USA}
\email{zanello@math.kth.se}

\thanks{${}^*$ The work for this paper was done while the first
author was sponsored by the National Security Agency under Grant
Number H98230-07-1-0036.\\\indent
2000 {\em Mathematics Subject Classification.} Primary: 13E10; Secondary: 13H10, 13D40.}

\begin{abstract}
We study a number of conditions on the Hilbert function of a level artinian algebra which imply the Weak Lefschetz Property (WLP).  Possibly the most important  open case is whether a codimension 3 SI-sequence forces the WLP for level algebras. In other words, does every codimension 3 Gorenstein algebra have the WLP?  We give some new partial answers to this old question: we prove an affirmative answer when the initial degree is 2, or when the Hilbert function is relatively small.  Then we give a complete answer to the question of what is the largest socle degree forcing the WLP.
\end{abstract}

\maketitle

\section{Introduction}

A very broad and fascinating problem in the study of standard graded algebras is to describe the algebraic, geometric and homological consequences forced on the algebras by conditions on the Hilbert function.  There is a vast literature on this topic.   This paper studies the behavior of Hilbert functions  of {\em artinian } algebras with respect to the Weak Lefschetz Property (WLP).  More specifically, we are interested in conditions on the Hilbert function  which {\em force} the WLP.  For  artinian algebras in general, this problem has already been solved  \cite{MZ1}.  Hence we refine it by focusing on level algebras, and we ask the following:

\medskip

\noindent {\bf Question:} {\em Let $\underline{h} = (1,h_1,h_2,\dots)$ be a Hilbert function that occurs for some artinian level algebra $A$.  What conditions on $\underline{h}$ guarantee that {\em every} level algebra with Hilbert function $\underline{h}$ has the WLP? More ambitiously, can we characterize these Hilbert functions? As an important special case, does every codimension 3 artinian Gorenstein algebra have WLP (i.e.\ is there an affirmative answer whenever $\underline{h}$ is a codimension 3  SI-sequence)?}

\medskip

Since it is not even known what Hilbert functions occur for level algebras in codimension $\geq 3$ (cf.\ \cite{GHMS, zanello2,weiss}), a full answer to this question will be very difficult to obtain.  This paper is intended to be a first step on the (possibly) long road toward its solution by giving a partial answer to the first part and to the last part.  It complements a result of Boij and the second author \cite{BZ}, who gave an example of a unimodal Hilbert function with the property that {\em no} level algebra  with this Hilbert function (which dictates that the algebra has type 2) has WLP.

Of course the last part of the question, the particular case of whether all codimension 3 artinian Gorenstein algebras have the WLP, has already been posed in the literature (see for instance \cite{HMNW, MMR3}).  It is known that in characteristic zero the WLP holds for a non-empty open subset of the codimension 3 Gorenstein algebras with fixed Hilbert function (e.g.\ \cite{harima, MN3}), and for {\em all} codimension 3 complete intersections \cite{HMNW} (it is false in positive characteristic -- cf.\ Remark \ref{char rmk}).  However, showing it for {\em all} codimension 3 Gorenstein algebras has proved very elusive so far, and we give some results in this direction.

An important prerequisite to answering these questions is to have a good understanding of as many Hilbert functions as possible for which non-WLP algebras exist.  As a first step, one can simplify the question by focusing only on the codimension, $r$, and socle type, $t$.  Are there pairs $(r,t)$ for which no non-WLP algebras exist?  As noted above, even taking codimension 3 and socle type 1 is open; this is a non-trivial question.

Indeed, much work has recently been done to prove the  existence of non-WLP level algebras for various $t$ and $r$. It is known that the WLP  holds for any standard graded algebra (level or not) when $r=2$ (see \cite{HMNW, MZ1}).  This topic was studied in \cite{GHMS}, where it was asked if at least for level algebras it might always hold in codimension 3 as well (\cite{GHMS},  Question 4.4).  The first counterexample was given by the second author \cite{zanello2}, who proved it for $t$ as small as 3.  This was extended  to $r=3, \  t=2$ in \cite{BZ}.  It has also  been shown that level algebras failing WLP exist for $r \geq 5$ and any $t$ (cf.\ \cite{stanley,BI,BL,boij, weiss}), as well as $r=3, \ t \geq 5$, and $r = 4, \ t \geq 3$ (cf.\ \cite{weiss}), and for $r=4, \  t=1$ (cf.\ \cite{ikeda, boij2}).

This leaves open only the  question of the existence of non-WLP level algebras for $r=3, \ t = 1$ (Gorenstein) and $t = 4$, and $r= 4, \ t = 2$. The latter two cases are quite easy to settle, just by variations of known constructions, as we will see in Example \ref{missing}.  The (major) case left open is therefore  only that of codimension 3 Gorenstein algebras, as we asked above.  
One goal of this paper is to begin the study of this latter problem.  

After some preparatory results and background material, this note is divided into two parts. In the first we study WLP for codimension 3 Gorenstein algebras, as indicated above. The first main result, Corollary \ref{init deg 2 has WLP}, will be that all such algebras enjoy the WLP if their initial degree is at most 2.  In the other main result of that section, Theorem \ref{small s}, we show the result for any socle degree, provided that the Hilbert function not be  too large.   One of the motivations of this section was to see to what extent the methods of \cite{MNZ3}, which the authors wrote together with Uwe Nagel, could be extended and modified to study WLP in codimension 3 rather than non-unimodality in codimension 4.

The second part of this paper brings the socle degree into the picture as well.  We  prove that the largest socle degree forcing all level algebras to enjoy the WLP is 2 if the codimension is $r= 3$, and  is 1 if $r\geq 4$.  One might wonder if this changes if we further restrict to codimension 3, type 2 level algebras, and indeed, we prove that the largest socle degree where all such algebras enjoy the WLP is 3, thus settling the closest case to codimension 3 Gorenstein algebras also with respect to the socle degree.


\section{Background and preparatory results}

Let $R = k[x_1, x_2,x_3]$ where $k$ is an infinite field.  For most of our results we need to assume that $k$ has characteristic zero, which we will see is indeed an essential hypothesis.

We first recall some standard terminology and notation.  Let $A$ be a standard graded artinian $k$-algebra, $A = \bigoplus_{i \geq 0} A_i$.  The {\em Hilbert function} of $A$ is the function $h_A$ defined by $h_A(i) = \dim_k A_i$.  The algebra $A$ has the {\em Weak Lefschetz Property (WLP)} if the homomorphism $(\times L) : A_i \rightarrow A_{i+1}$ induced by multiplication by a general linear form $L$ has maximal rank for all $i$.  It has the {\em Strong Lefschetz Property (SLP)} if the homomorphism $(\times L^s) : A_i \rightarrow A_{i+s}$ has maximal rank for all $i$ and all $s$. 
We say that $A$ is {\em level of type $t$} if the socle of $A$ is of dimension $t$ and is concentrated in one degree (namely the last degree in which $A$ is non-zero).  Furthermore,  $A$ is {\em Gorenstein} if and only if it is level of type 1.

Three basic results studying the behavior of Hilbert functions are those of Macaulay, Gotzmann and Green, for which we need a little combinatorial notation first.

\begin{definition} Let $n$ and $i$ be positive integers. The {\em i-binomial expansion of n} is 
\[
n_{(i)} = \binom{n_i}{ i}+\binom{n_{i-1}}{ i-1}+...+\binom{n_j}{ j},
\]
 where $n_i>n_{i-1}>...>n_j\geq j\geq 1$. Such an expansion always exists and it is unique (see, e.g.,  \cite{BH}, Lemma 4.2.6).\\\indent
Following \cite{BG}, we define, for any integers $a$ and $b$,
$$
(n_{(i)})_{a}^{b}=\binom{n_i+b}{ i+a}+\binom{n_{i-1}+b}{ i-1+a}+...+\binom{n_j+b}{ j+a},
$$
where we set $\binom{m}{ q}=0$ whenever $m<q$ or $q<0$.
\end{definition}

\begin{theorem} Let $L \in A$ be a general linear form. Denote by $h_d$  the degree $d$ entry of the Hilbert function of $A$ and by $h_d^{'}$ the degree $d$ entry of the Hilbert function of $A/L A$. Then:\\\indent
i) {\em (Macaulay)} $$h_{d+1}\leq ((h_d)_{(d)})_{1}^{1}.$$\indent
ii) {\em (Gotzmann)} 
If $h_{d+1} = ((h_d)_{(d)})^1_1$ and $I$ is generated in degrees $\leq d$ then 
\[
h_{d+s} = ((h_d)_{(d)})^s_s \hbox{\hspace{.3cm} for all $s \geq 1$}.
\]
\indent iii) {\em (Green)} $$h_d^{'}\leq ((h_d)_{(d)})_{0}^{-1}.$$
\end{theorem}

\begin{proof}
 i) See \cite{BH}, Theorem 4.2.10.\\
\indent ii) See \cite{BH}, Theorem 4.3.3, or \cite{Go}. \\
\indent iii) See \cite{Gr}, Theorem 1. 
\end{proof}

A sequence of non-negative integers $\underline{h} = (1,r,h_2, \dots , h_d ,\dots )$ is said to be an {\em O-sequence} if it satisfies Macaulay's bound for all $d$.  We remark that Macaulay also showed that any O-sequence is actually the Hilbert function of some standard graded algebra, so the O-sequences are precisely the Hilbert functions of standard graded algebras.    When $A$ is artinian and Gorenstein, its Hilbert function is a symmetric O-sequence.  Such a sequence is a {\em Stanley-Iarrobino (SI)-sequence} if in addition, the first half is  {\em differentiable}, i.e.\ its first difference is also an O-sequence.  Such sequences are necessarily unimodal.

\begin{lemma}[\cite{MMN}]  \label{gen wlp}
Let $R/I$ be an artinian standard graded algebra and let $L$ be a general linear form.  Consider the homomorphisms $\phi_d : (R/I)_d \rightarrow (R/I)_{d+1}$ defined by multiplication by $L$, for $d \geq 0$.  Note that $(R/I)_d$ and $(R/I)_{d+1}$ are finite-dimensional vector spaces.

\begin{itemize}
\item[(a)] If $\phi_{d_0}$ is surjective for some $d_0$ then $\phi_d$ is surjective for all $d \geq d_0$.

\item[(b)] If $R/I$ is level and $\phi_{d_0}$ is injective for some $d_0 $,  then $\phi_d$ is injective for all $d \leq d_0$.

\item[(c)] In particular, if $R/I$ is level and $\dim (R/I)_{d_0} = \dim (R/I)_{d_0+1}$ for some $d_0$ then $R/I$ has WLP if and only if $\phi_{d_0}$ is injective (and hence is an isomorphism).
\end{itemize}
\end{lemma}

\begin{remark}\label{rrr}
Lemma \ref{gen wlp} implies that for Gorenstein algebras, there is always exactly one degree that needs to be checked, and it can be chosen so that only injectivity (resp.\ surjectivity) has to be checked.  Indeed, the only missing ingredient is that since $R/I$ is self-dual up to twist, injectivity in the ``first half'' is equivalent to surjectivity in the ``second half.''  In Section \ref{gor alg sect}, we will refer to  this with the phrase ``by duality."
\end{remark}

\begin{proposition}[\cite{MNZ3}] \label{prop-key}
Assume that the field $k$ has characteristic zero.
Let $R =$ \linebreak $ k[x,y,z]$ and let $J = (F,G_1,G_2) \subset R$ be a homogeneous ideal with three minimal generators, where $\deg F = a \geq 2$ and $\deg G_1 = \deg G_2 = b \geq  a$.  Let $L  \in R$ be general linear form.  Then $\dim [R/(J,L)]_b = a-1$ if and only if $F,G_1,G_2$ have a GCD of degree $a-1$.  Otherwise $\dim [R/(J,L)]_b = a-2$.
\end{proposition}

\begin{proof}
We only remark that this result is stated in \cite{MNZ3} in a slightly different way, but in the proof it immediately passes to this setting.
\end{proof}

\begin{lemma} \label{gcd observations}
Let $R/I$ be an artinian graded algebra, and let $a = \min \{ t \  | \ I_t \neq 0 \}$ be the initial degree of $I$.  Let $F \in I$ be a form of degree $a$ and $L$ a general linear form.  Let $d \geq a$ be an integer.

\begin{enumerate}
\item  If $I_d$ has a GCD of degree $a-\epsilon$ then $\dim [R/(I,L)]_d \geq a-\epsilon$.

\item Suppose that $\dim [R/(I,L)]_d > a-\delta$.   Suppose that $I$ has some minimal generating set that contains $F$ together with $\delta$ generators of degree $d$ (not including $F$ if $a = d$), and contains no additional minimal generators of degree $<d$. Then  the multiplication by a general linear form, 
\[
(\times L) : (R/I)_{d-1} \rightarrow (R/I)_d
\]
fails to be injective.  

\end{enumerate}
\end{lemma}

\begin{proof}
For (1), let $G$ be the GCD of $I_d$ and let $\bar G$ be its image in $R/(L)$.  Then  the elements of $(I,L)$ in any degree $t \leq d$, viewed in $R/(L)$, are all of the form $\bar G f$ where $f$ is a form in $[R/(L)]_{t-(a-\epsilon)} \cong k[x,y]_{t-(a-\epsilon)}$.  Note that $\bar F$ is one such element, with $t = a$.  Since $\dim [R/(L)]_d = d+1$ and since there are $(d-a+\epsilon +1)$ independent forms of degree $d-a+\epsilon$ in $k[x,y]$, we see that the maximum number of independent elements of $(I,L)/(L)$ in $R/(L)$ is $d-a+\epsilon+1$, so $\dim[R/(I,L)]_d \geq d+1 - (d-a+\epsilon+1) = a-\epsilon$ as claimed.

For (2), note that $I$ may have more than $\delta$ minimal generators of degree $d$; we just require the condition on the Hilbert function. Let $F_1,\dots,F_\delta$ be the indicated minimal generators.  Then the hypotheses say that $\bar F, \bar F_1,\dots,\bar F_\delta$ are linearly dependent modulo $L$.  This means that there is some form, $A$, of degree $d-1$ such that $AL + a_1F_1 + \dots + a_\delta F_\delta = 0$.  Since the $F_i$'s are homogeneous, for degree reasons clearly $A \notin \langle F_1,\dots,F_\delta \rangle$.  If $A \in (F)$ then one of the $F_i$ is redundant, contradicting their choice.  This means that $AL \in I$ but $A \notin I$, so $A$ is a nonzero element in the kernel of the multiplication.
\end{proof}

\begin{remark}
Lemma \ref{gcd observations} explains {\em why} Example 7 of  \cite{zanello2} works.
\end{remark}

\begin{remark}  We collect the following easy facts.

\begin{enumerate}
\item Let $I$ be  any homogeneous ideal  and let $F$ be any form of degree $d$. There is an exact sequence:





\begin{equation} \label {usual exact seq}
0 \rightarrow R/(I:F)(-d) \stackrel{\times F}{\longrightarrow} R/I \rightarrow R/(I,F) \rightarrow 0.
\end{equation}


\medskip

\item If $R/I$ is Gorenstein then we have the following well-known facts:

\medskip

\begin{itemize}

\item  If $F\notin I$, then $R/(I:F)$ is Gorenstein, of socle degree $e - d$.

\medskip

\item For any form $G\notin I$ of degree $g \geq 1$,  $h_{R/(I:G)} (d-g) \leq h_{R/(I:G)}(d-g+1).$
Indeed, codimension 3 Gorenstein Hilbert functions are always SI-sequences, hence unimodal (cf.\ \cite{stanley} and the recent elementary proof in \cite{zanello1}), so it is enough to check that $d-g+1 \leq \frac{e-g}{2}$, which is an easy calculation (since $g \geq 1$).  

\end{itemize}
\end{enumerate}
\end{remark}

We now need some preliminary lemmas. A very useful rephrasing of the Weak Lefschetz Property in terms of the exact sequence (\ref{usual exact seq}) is the following:

\begin{lemma}\label{lll}
An artinian algebra $R/I$ with Hilbert function $(1,h_1,h_2,...,h_e)$ enjoys the WLP if and only if, for any general form linear $L$ and for all indices $i$, we have
$$h_{R/(I,L)}(i)=\max \{h_i-h_{i-1} ,0\}.$$
\end{lemma}

Notice, also, that in the level (or Gorenstein) case, by Lemma \ref{gen wlp}, in order to check whether WLP holds, there always exists one degree such that it suffices to check the value of $h_{R/(I,L)}(i)$ only in that spot.

\begin{lemma} \label{force surjective}
If $R/I$ is any standard graded algebra and $L$ is a linear form, then $(R/I)_t \stackrel{\times L}{\longrightarrow} (R/I)_{t+1}$ is surjective if and only if $(R/(I,L))_{t+1} = 0$.
\end{lemma}

\begin{proof}
The exact sequence (\ref{usual exact seq}) applied to $F = L$ immediately gives the result.
\end{proof}

\begin{lemma}\label{force gcd}
If $h_{R/(I,L)}(t) = h_{R/(I,L)}(t+1) = k$ then $I_{t+1}$ has a GCD of degree $k$.
\end{lemma}

\begin{proof}
This simple but powerful tool was used in \cite{MNZ3}.  The point is that Davis's theorem forces a GCD in $(I,L)_{t+1}$, which then lifts to $I_{t+1}$ (cf.\ \cite{davis, BGM}).
\end{proof}

\begin{lemma} \label{no GCD}
Let $R/I$ be an artinian Gorenstein algebra of socle degree $e$.  If $d > \frac{e}{2}$ then $I$ does not have a GCD of any degree $r \geq 1$ occurring in degree $d$.
\end{lemma}

\begin{proof}
Suppose otherwise, and let $F$ be such a GCD.  Now we apply (\ref{usual exact seq}) with the GCD, $F$, playing the role of the homogeneous form.  Note that $(R/(I,F))_t = (R/(F))_t$ for $t \leq d$.  Note also that $R/(I:F)$ is Gorenstein of socle degree $e-r$.  Finally, since $I_t$ has to have (at least) $F$ as a common divisor in all degrees $\leq d$, without loss of generality we will assume that $d = \lfloor \frac{e+2}{2} \rfloor$.

Then we have
\[
\begin{array}{rcll}
h_{R/I}(d-1) & \geq & h_{R/I}(d) & \hbox{by definition of $d$} \\
h_{R/(I:F)(-r)}(d-1) & \leq & h_{R/(I:F)(-r)}(d)  \\
h_{R/(I,F)}(d-1) & < & h_{R/(I,F)}(d) & \hbox{since $(I,F) = (F)$ in this range}
\end{array}
\]
where the second inequality follows because of the (revised) definition of $d$ and because all Gorenstein Hilbert functions are unimodal in codimension 3.
But then
\[
h_{R/I}(d-1) = h_{R/(I:F)(-r)}(d-1) + h_{R/(I,F)}(d-1) < h_{R/(I:F)(-r)}(d) + h_{R/(I,F)}(d) = h_{R/I}(d)
\]
is a contradiction. 
\end{proof}

\begin{remark}
It is an open question whether all Gorenstein Hilbert functions in codimension 4 are unimodal.  This was shown in \cite{MNZ3} for $h_4 \leq 33$.  If it is true that all such Hilbert functions are unimodal, then Lemma \ref{no GCD} holds in codimension 4 as well.
\end{remark}


\section{Gorenstein algebras of codimension 3} \label{gor alg sect}

We begin with a useful result connecting WLP with GCD's of components of ideals.

\begin{lemma} \label{GCD implies WLP}
Let $R/I$ be an artinian Gorenstein algebra of socle degree $e$.  Set $d = \lfloor \frac{e-1}{2} \rfloor$. Let $a = \min \{ t \ | \  I_t \neq 0 \}$ (the initial degree of $I$).  If $I_{d+1}$ has a GCD, $F$, of degree $a- 1$ then $e$ is even and  $R/I$ has WLP.
\end{lemma}

\begin{proof}
If $e$ is odd then $d+1 > \frac{e}{2}$, so Lemma \ref{no GCD} shows that no such algebras exist.  Hence $e$ must be even.

 In order to show that $R/I$ has WLP, it suffices to check that the multiplication by a general linear form from degree $d$ to degree $d+1$ is injective, by Lemma \ref{gen wlp} and duality.

\bigskip

We make the following observations:

\begin{enumerate}
\item \label{F} For any $t \leq d+1$, $(I,F)_t = (F)_t$.

\item \label{init deg} $(I:F)$ has initial degree 1.  Hence $R/(I:F)$ is isomorphic to a codimension two Gorenstein algebra (necessarily a complete intersection).
\end{enumerate}

Now, the exact sequence (\ref{usual exact seq}) gives rise to the following diagram (after taking into account observation (\ref{F}) and the previous claim):
\[
\begin{array}{ccccccccccccc}
0 & \rightarrow & (R/(I:F))_{d-(a-1)} & \rightarrow & (R/I)_d & \rightarrow & (R/(F))_d & \rightarrow & 0 \\
&& \downarrow && \downarrow && \downarrow \\
0 & \rightarrow & (R/(I:F))_{d-(a-1)+1} & \rightarrow & (R/I)_{d+1} & \rightarrow & (R/(F))_{d+1} & \rightarrow & 0
\end{array}
\]
where the vertical arrows are multiplication by a general linear form.  The leftmost vertical map is injective by (\ref{init deg}), and the rightmost vertical map is clearly injective.  Thus the middle map is injective, and so $R/I$ has WLP.
\end{proof}

\begin{corollary} \label{init deg 2 has WLP}
If $R/I$ is Gorenstein and $I$ has initial degree 2, then $R/I$ has WLP.
\end{corollary}

\begin{proof}
We consider the possibilities for $h_2 = h_{R/I}(2)$.  

\bigskip

\noindent \underline{\it Case 1}: $h_2 = 3$.  Since the Hilbert function of $R/I$ is an SI-sequence, it  is of the form ${\tt 1,3,3,3, \dots,3,3,1}$.  Suppose first  that the socle degree is $\geq 4$.  Then in particular, $h_{R/I}(3) = 3$.  By Green's theorem, $h_{R/(I,L)}(3) = 0$, so by Lemma \ref{force surjective}, we have a surjectivity $(R/I)_t \rightarrow (R/I)_{t+1}$ for all $t \geq 2$.  In particular, since at least $h_{R/I}(3) = 3$, by duality we conclude that $R/I$ has WLP.

It remains to prove WLP for a Gorenstein algebra with Hilbert function ${\tt 1,3,3,1}$.  Let $F$, $G_1$ and $G_2$ all be minimal generators of $I$ of degree 2, so $a = b = 2$ in Proposition \ref{prop-key}.  Suppose that $R/I$ fails to have WLP.  Then the multiplication by a general linear form from degree 1 to degree 2 fails to be surjective.  By Lemma \ref{force surjective}, this means that $h_{R/(I,L)}(2) = 1 = 2-1$.  Hence by Proposition \ref{prop-key},  $I_2$ has a GCD of degree $1 = 2-1$.  Then by Lemma \ref{GCD implies WLP}, $R/I$ has WLP. We only remark that, as pointed out to us by the referee,  WLP for Gorenstein algebras with $h_2 = 3$ can also be proved, still assuming that the characteristic of the base field be zero, using an argument involving the Hessian of an inverse system form (substantially due to \cite{GN} and \cite{Wa}, and also employed in \cite{HMNW}, Example 4.3).  See also Remark \ref{hessian}.
\bigskip

\noindent \underline{\it Case 2}: $h_2 = 4$.  First suppose that the Hilbert function has of one of the following forms:
\[
\begin{array}{l}
{\tt 1,3,4,3,1} \\
{\tt 1,3,4,4,\dots,4,3,1} \\
{\tt 1,3,4,5,5, \dots ,5,4,3,1}.
\end{array}
\]
Using an argument almost identical to the one for Case 1, we get that multiplication by a general linear form is surjective from degree 2 to degree 3 in the first two cases, and from degree 3 to degree 4 in the third case.  All of these are enough to force WLP.

It remains to consider the case ${\tt 1,3,4,5,6,\dots}$.  Now there are two possibilities: 
\[
\begin{array}{l}
{\tt 1,3,4,5,6,\dots,t-1,t,t-1,\dots} \\
{\tt 1,3,4,5,6,\dots,t-1,t,t, \dots,t,t-1,\dots}
\end{array}
\]
In the first of these cases, the second $t-1$ occurs in degree $t-1$, and Green's theorem together with Lemma \ref{force surjective} guarantees that multiplication by a general linear form from degree $t-2$ to $t-1$ is surjective.  By Lemma \ref{gen wlp} and duality, this implies that $R/I$ has WLP.  In the second of these cases, a similar argument takes care of the case where there are at least three $t$'s.  So we have to check the case ${\tt 1,3,4,5,6,\dots,t-1,t,t,t-1,\dots}$.  Note that the second $t$ occurs in degree $t-1$.  Now Green's theorem gives that $h_{R/(I,L)}(t-1) \leq 1.$  If it is equal to 0, then again we have WLP.  So without loss of generality assume that it is 1.  But also applying Green's theorem to degree $t-2$, we obtain that $h_{R/(I,L)}(t-2) = 1$.  Then Lemma \ref{force gcd} gives that $I_{t-1}$ has a GCD of degree 1, so by Lemma \ref{GCD implies WLP}, $R/I$ has WLP.

\bigskip

\noindent \underline{\it Case 3}: $h_2 = 5$.  Now the form of the Hilbert function is essentially the following:
\[
{\tt 1, 3, 5, \hbox{(grow by 2)}, \hbox{(grow by 1)}, \hbox{(flat)}, \dots}
\]
where any of these ranges could be empty.  All of the specific subcases are dealt with using the same ideas as above, and we omit the details except for one that has a slight twist.  Suppose that the Hilbert function is of the form
\[
{\tt 1,3,5,7,\dots,2t-1, 2t+1, 2t+1, 2t-1,\dots}
\]
where the first $2t+1$ occurs in degree $t$.  We then note that $I$ has a generator in degree $a=2$ and two in degree $b = t+1$, so by Lemma \ref{force surjective} and Proposition \ref{prop-key},  $R/I$ fails to have WLP if and only if $I_{t+1}$ has a GCD of degree 1, and Lemma \ref{GCD implies WLP} gives the result.
\end{proof}

\begin{remark} \label{char rmk}
The proof for the case ${\tt 1,3,3,1}$ used Proposition \ref{prop-key} in a crucial way (although,  as we said above, a Hessian argument can also be used). An important hypothesis, in either case, is that $k$ has characteristic zero.  It was pointed out to us by Uwe Nagel that  in fact in characteristic 2 this Hilbert function does {\em not} necessarily have WLP: the complete intersection $(x_1^2,x_2^2,x_3^2)$ is a counterexample (cf.\ \cite{HMNW}).  This answers a question raised  at the beginning of section 2 of \cite{MNZ3}, whether Proposition \ref{prop-key} has a characteristic-free proof.  (In fact, a counterexample can be found in any characteristic, using a complete intersection of forms of the same degree, via the same approach.)
\end{remark}

\begin{corollary}
Let $R/I$ be an artinian Gorenstein algebra such that $I$ has minimal generators $F, F_1,F_2$ of degrees $2,b,b$ respectively ($b \geq 2$),  $F,F_1,F_2$ have  a GCD in degree $b$, and all other generators of $I$ have degree $\geq b$.  Then $R/I$ cannot have WLP.  
\end{corollary}

\begin{proof}
Since the initial degree of $I$ is $2$, we see that any GCD would have to have degree 1.  By Lemma \ref{no GCD}, the socle degree must be $\geq 3$ and $h_{R/I}(b-1) \leq h_{R/I}(b)$ (since otherwise $I$ has a GCD in degree $> \frac{e}{2}$).  So WLP would mean in particular that we need injectivity from degree $b-1$ to degree $b$.  By Lemma \ref{gcd observations} (1), $\dim [R/(I,L)]_b \geq 2-1=1$.  Then by Lemma \ref{gcd observations} (2) (taking $\delta = 2$) the required injectivity fails.
\end{proof}

Corollary \ref{init deg 2 has WLP} allows us to extend Lemma \ref{GCD implies WLP}, lowering by 1 the degree of the GCD that forces WLP:

\begin{corollary}
Let $R/I$ be an artinian Gorenstein algebra of socle degree $e$.  Set $d = \lfloor \frac{e-1}{2} \rfloor$. Let $a = \min \{  t \  | \  I_t \neq 0 \}$ (the initial degree of $I$).  If $I_{d+1}$ has a GCD of degree $a- 2$ then $e$ is even and  $R/I$ has WLP.
\end{corollary}

\begin{proof}
The proof is almost identical to that of Lemma \ref{GCD implies WLP}.  The only difference is that now with the GCD, $F$, of degree $a-2$, we get that $I:F$ has initial degree 2 rather than 1.  But then with the same argument, also invoking Corollary \ref{init deg 2 has WLP}, we obtain the result.
\end{proof}

\begin{theorem} \label{small s}
Let $R/I$ be a Gorenstein artinian algebra with socle degree $e$ and Hilbert function $h_i = h_{R/I}(i)$.  Assume that there is some  integer $s$ such that
\begin{itemize}
\item $3 \leq s \leq \frac{e}{2}-1$;
\item $h_s \leq 3s-1$;
\end{itemize}
Then $R/I$ has WLP.
\end{theorem}

\begin{proof}
Note that the $s$-binomial expansion of $3s-1$ is
\[
3s-1 = \binom{s+1}{s} + \binom{s}{s-1} + \binom{s-2}{s-2} + \dots + \binom{1}{1}.
\]
Then the condition that $h_s \leq 3s-1$ implies, by Green's theorem, that $h_{R/(I,L)}(i) \leq 2$ for all $i \geq s$.  

Suppose first that $e$ is odd, and set $d = \frac{e-1}{2}$.  Then $s \leq d-1$, and $h_{d-1} = h_{d+2} \leq h_d = h_{d+1}$.  The failure of WLP would imply that $h_{R/(I,L)}(d+1) \geq 1$.  Since $h_{R/(I,L)}(d-1) \leq 2$,  $h_{R/(I,L)}(d)$ is equal to either $h_{R/(I,L)}(d-1)$ or $h_{R/(I,L)}(d+1)$ (or both).  In the latter case, $I$ has a GCD in degree $d+1 > \frac{e}{2}$, which is impossible by Lemma \ref{no GCD}.  

So without loss of generality  suppose that $h_{R/(I,L)}(d-1) = h_{R/(I,L)}(d) > h_{R/(I,L)}(d+1)$.  This can only happen if these values are 2, 2 and 1 respectively.  Hence there is a GCD, $Q$, of degree 2 in $I_d$.  Reducing modulo a general linear form $L$, we observe that $(I,L)_t = (Q,L)_t$ for $t = d-1$ and $d$ since one inclusion is clear and they have the same Hilbert function in those degrees.  

We now consider the other relevant Hilbert functions.  For clarity we separate the steps.

\begin{enumerate}
\item $h_{R/(I,L)}(d+1) = 0$ since the only other possibility is that it equals 1, but then $I$ has a GCD of degree 1 in degree $d+1$, violating Lemma \ref{no GCD}.  

\item From (\ref{usual exact seq}) and the values obtained above we get
\[
\begin{array}{rcl}
h_{R/(I:L)}(d-2) & = & h_{d-1} - 2 \\
h_{R/(I:L)}(d-1) & = & h_d - 2 \\
h_{R/(I:L)}(d) & = & h_{d+1} -1 \\
h_{R/(I:L)}(d+1) & = & h_{d+2} 
\end{array}
\]

\item By the symmetry of $h_{R/(I:L)}$, we have that $h_{R/(I:L)} (d-1) = h_{R/(I:L)} (d+1)$.  Hence from the equalities above and the symmetry of $h_{R/I}$, 
\[
h_{d-1} = h_{d+2} = h_{R/(I:L)} (d+1) = h_{R/(I:L)} (d-1)  = h_d -2 .
\]
\end{enumerate}

This last equality shows that the Hilbert function of $R/I$ grows by 2 from degree $d-1$ to degree $d$.  The binomial expansion above, together with Macaulay's theorem, imply that this growth is maximal.   By Gotzmann's theorem, this implies that $I_d$ has a GCD, $Q$, of degree 2 (as we saw above) and furthermore that if we set $J = \langle I_{d-1} \rangle$ to be the ideal generated by the homogeneous component of degree $d-1$, the Hilbert function of  $R/J$ grows by two in all subsequent degrees.  Since $h_{d+1} = h_d$, this means that $I$ has exactly two new generators in degree $d+1$, say $F_1$ and $F_2$.

Now consider the ideal $(Q,F_1,F_2)$.  We have seen that $(I,L)_d = (Q,L)_d$.  Hence $(Q,F_1,F_2, L)_{d+1} = (I,L)_{d+1}$, and so $1 = h_{R/(I,L)}(d+1) = h_{R/(Q,F_1,F_2,L)}(d+1)$.  But this is exactly the situation of Proposition \ref{prop-key}, and it implies that $Q$, $F_1$ and $F_2$ have a GCD of degree 1.  This means that $I_{d+1}$ has a GCD, which violates Lemma \ref{no GCD}.  This completes the case that $e$ is odd.

Now let the socle degree $e$ be even. It suffices to show that the multiplication by a general linear form $L$ is injective between degrees $d=\frac{e}{2}-1$ and $d+1=\frac{e}{2}$. If $h_d=h_{d+1}$, then by symmetry $h_d=h_{d+1}=h_{d+2}$, and the WLP follows from a result of Iarrobino and Kanev (\cite{IK} Theorem 5.77): they show that, if the Hilbert function of a codimension 3 Gorenstein algebra $R/I$ has three consecutive entries $a,a,a$, then there is a unique zero-dimensional subscheme of $\mathbb P^2$ of degree $a$ whose ideal is equal to $I$ in the three degrees where $I$ has dimension $a$.  But clearly then $R/I$ has depth 1, so multiplication by a general linear form in those degrees is injective and we have the WLP for $R/I$.

So we may suppose that $h_{d+1}>h_d$. Similarly to what we have observed above, by Macaulay's theorem  we have that either  $h_{d+1}=h_d+1$ or  $h_{d+1}=h_d+2$. Likewise, since $d\geq s$, $h_{R/(I,L)}(d+1)\leq h_{R/(I,L)}(d)  \leq 2$.

If $h_{R/(I,L)}(d+1) =1$, then we clearly must have  $h_{d+1}=h_d+1$, and the WLP follows. Thus it remains to consider when $h_{R/(I,L)}(d+1)= h_{R/(I,L)}(d)  = 2$. In this case, we have a degree 2 GCD, say $Q$, for $I_{d+1}$. An argument, involving the unimodality of  $R/(I:Q)$, entirely similar to the one we gave for the previous case implies that $h_{d+1} =h_d+2$. But since $h_{R/(I,L)}(d+1)  = 2$, we have that the multiplication by $L$ between $(R/I)_{d}$ and $(R/I)_{d+1}$ is injective, that is that $R/I$ has the WLP, as desired.
\end{proof}

\begin{corollary}
Let $R/I$ be a Gorenstein artinian algebra with $h_3 \leq 8$.  Assume that the Hilbert function of $R/I$ is not ${\tt 1,3,6,8,8,6,3,1}$ or ${\tt 1,3,6,6,3,1}$.  Then $R/I$ has WLP.
\end{corollary}

\begin{proof}
Note that if $h_2 < 6$ we already know the result from Corollary \ref{init deg 2 has WLP}.  So without loss of generality we assume that $h_2 = 6$.  

The condition $h_3 \leq 8$ implies (by Macaulay's theorem) that $h_s \leq 2s+2$ for all $s \geq 3$.  Theorem \ref{small s} assumes that $e \geq 8$.   Hence we only have to take care of the cases involving small socle degree. The smallest possibility is $e=6$.  For convenience, we will summarize the numerical information obtained from the exactness of (\ref{usual exact seq}) in a table.  We set $L$ to be a general linear form.

\bigskip

\noindent \underline{\it Case 1}: ${\tt 1,3,6,8,6,3,1}$.

Green's theorem gives that $h_{R/(I,L)}(3) \leq 2$. It follows that the only possible values of the corresponding Hilbert functions are 

\begin{center}
\begin{tabular}{c|cccccccccccccccccccccccccc}
deg & 0 & 1 & 2 & 3 & 4 & 5 & 6 \\ \hline
$h_{R/I}$ & 1 & 3 & 6 & 8 & 6 & 3 & 1 \\
$h_{R/(I:L)(-1)}$ & & 1 & 3 & 6 & 6 &3 & 1  \\ \hline
$h_{R/(I,L)}$ & 1 & 2 & 3 & 2 & 0 & 0 & 0
\end{tabular}
\end{center}

Hence this case follows from Lemma \ref{lll}.

\bigskip

\noindent \underline{\it Case 2}: ${\tt 1,3,6,7,6,3,1}$.

\begin{center}
\begin{tabular}{c|cccccccccccccccccccccccc}
deg & 0 & 1 & 2 & 3 & 4 & 5 & 6 \\ \hline
$h_{R/I}$ & 1 & 3 & 6 & 7 & 6 & 3 & 1 \\
$h_{R/(I:L)(-1)}$ & & 1 & 3 && &3 & 1  \\ \hline
$h_{R/(I,L)}$ & 1 & 2 & 3 
\end{tabular}
\end{center}

If $R/I$ fails WLP then $h_{R/(I,L)}(4) \geq 1$.  But Green's theorem applied to degree 3 implies that $h_{R/(I,L)}(3) \leq 2$.  Hence there are two possibilities.  If $h_{R/(I,L)}(3) = h_{R/(I,L)}(4) = 1$ then $I$ has a GCD of degree 1 in degree 4, and we conclude with Lemma \ref{GCD implies WLP} (or just observe directly that the multiplication from degree 2 to degree 3 is injective, which is enough).  If $h_{R/(I,L)}(3) = 2$ then the three generators of $I$ in degree 3 fail to be independent modulo $L$, so Proposition \ref{prop-key} applies.  It is impossible for $h_{R/(I,L)}(3)$ to equal 0.

\bigskip

\noindent \underline{\it Case 3}: ${\tt 1,3,6,6,6,3,1}$.
The three consecutive 6's imply WLP by \cite{IK}, Theorem 5.77.

\bigskip

\noindent \underline{\it Case 4}: ${\tt 1,3,6,6,6,6,3,1}$.
This case is immediate using these methods.

\bigskip

\noindent \underline{\it Case 5}: ${\tt 1,3,6,7,7,6,3,1}$.
This case is immediate using these methods.

\bigskip

We remark that for each of the two missing  cases, the considerations above leave only  one possibility: 

\begin{center}
\begin{tabular}{c|ccccccccccccc|cccccccccccccccccc}
deg & 0 & 1 & 2 & 3 & 4 & 5 & 6 & 7 &&&&& deg & 0 & 1 & 2 & 3 & 4 & 5 \\ \cline{1-9} \cline{14-20}
$h_{R/I}$ & 1 & 3 & 6 & 8 & 8 & 6 & 3 & 1 &&&&& $h_{R/I}$ & 1 & 3 & 6 & 6 & 3 & 1 \\
$h_{R/(I:L)(-1)}$ & & 1 & 3 & 6 & 7 & 6 &3 & 1 &&&&& $h_{R/(I:L)}(-1)$ && 1 & 3 & 5 & 3 & 1   \\ \cline{1-9} \cline{14-20}
$h_{R/(I,L)}$ & 1 & 2 & 3 & 2 & 1 &&& &&&&& $h_{R/(I,L)}$ & 1 & 2 & 3 & 1
\end{tabular}
\end{center}

In the first of these, if $R/I$ is a complete intersection (of type (3,3,4)), then WLP is known by \cite{HMNW}. 
\end{proof}

As mentioned in the introduction, the most natural (and most important) question  at this point is the following: 

\begin{question}
{\em Do all codimension 3 Gorenstein algebras possess the WLP?}
\end{question}


\section{Small socle degree}

Let us now turn our attention to the problem of determining the largest socle degree forcing the WLP for all level algebras of any given codimension, as well as for some interesting specific cases, such as those of level algebras of codimension 3 and type 2 and 3. We refer to \cite{Ge,IK} for an introduction to the theory of Macaulay's inverse systems, which will be needed in this portion of the paper.

We define $e(r)$ as the largest socle degree $e$ such that all level algebras of codimension $r$ and socle degree $\leq e$ enjoy the WLP (putting $e(r)=+\infty $ if such integer does not exist). Also, set $e_t(r)$ to be the analogous value when we restrict to type $t$.

We begin with the following construction, which proves the existence of a type 2, codimension 3 level algebra of socle degree 4. It is motivated by an inspiring  example of Brenner-Kaid. 

\begin{lemma} \label{brenner} Let $A=k[x_1,x_2,x_3]/I$ be the codimension 3 level algebra corresponding to the inverse system module  $M=\langle y_1^2(y_2^2+y_3^2),y_2^2(y_1^2+y_3^2)\rangle \subset k[y_1,y_2,y_3]$. Then $A$ has Hilbert function $(1,3,6,6,2)$ and  fails to have the WLP. In particular, $e_2(3)\leq 3$.
\end{lemma}

\begin{proof}
Brenner and Kaid  (\cite{BK}, Example 3.1) proved that the artinian algebra
\[
k[x_1,x_2,x_2]/(x_1^3,x_2^3,x_3^3,x_1x_2x_3),
\]
which has  Hilbert function $(1,3,6,6,3)$, fails to have the WLP between degree 2 and 3. It is easy to see that this is a level algebra of type 3, for instance by computing its inverse system module  (it is immediate with CoCoA \cite{CoCoA}), which is $M'= \langle y_1^2y_2^2,y_2^2y_3^2, y_1^2y_3^2 \rangle \subset k[y_1,y_2,y_3]$.

Now, by computing the first partial derivatives of the generators of $M'$, we see that (as a $k$-vector space) 
\[
M'_3= \langle y_1^2y_2,y_1^2y_3,y_1y_2^2,y_1y_3^2,y_2^2y_3,y_2y_3^2 \rangle.
\]
Since both $M$ and $M'$ are level algebras, it is enough to prove that $M_3=M'_3$ in order to conclude that the algebra $A$ of the statement has Hilbert function $(1,3,6,6,2)$ and also fails to have the WLP (from degree 2 to 3), because such equality on the inverse systems implies that the two corresponding ideals also coincide in degrees $\leq 3$. But that is a standard computation of linear algebra, so will be omitted.
\end{proof}

The next construction provides the existence of socle degree 2, level algebras of codimension $\geq 4$ without the WLP.

\begin{lemma} \label{222} Let $r\geq 4$. The level algebras quotients of $k[x_1,...,x_r]$, whose inverse system module is generated by $M= \langle y_1^2,y_1y_2,y_2^2,y_3y_4,y_5^2,...,y_r^2 \rangle$, all fail to have the WLP. In particular, $e(r)\leq 1$ for $r\geq 4$.
\end{lemma}

\begin{proof} Note first of all that the algebras of the statement have Hilbert function $(1,r,r)$, since $M$ is generated by $r$ monomials of degree 2 and by differentiating them we obtain all the variables. 

Also, it is easy to show that (unlike the case of arbitrary ideals), if an ideal $I$  in some degree $i$ is spanned, as a $k$-vector space, by monomials, then the subgroup $M_i$ of the inverse system $M$ of $I$ is also monomial, and furthermore its monomials are the same as those whose classes span $(R/I)_i$ (of course after renaming the variables). Therefore, in degrees where an ideal is monomial, proving the WLP can be done on the inverse system.

Thus, it is enough to show that, for any linear form $L=a_1y_1+...+a_ry_r$, the multiplication by $L$ between $M_1= \langle y_1,y_2,...,y_r \rangle$ and $M_2= \langle y_1^2, y_1y_2, y_2^2, y_3y_4, y_5^2,...,y_r^2 \rangle$ is not a bijective map (thinking of the generators of $M_1$ and $M_2$ as those of $(R/I)_1$ and $(R/I)_2$). Suppose it is. Then it is also injective, and since a standard computation shows  that $a_3y_3-a_4y_4$ is in the kernel, we must have $a_3=a_4=0$. But then it easily follows that $L\cdot y_3=L\cdot y_4=0$, a contradiction.
\end{proof}

We now have a key lemma, whose argument relies on those of the previous section:

\begin{lemma}\label{133}
All level algebras whose Hilbert functions start $(1,3,3)$ enjoy the WLP.
\end{lemma}

\begin{proof} 
As before, we take $F,G_1,G_2$ to be minimal generators of degree 2. By Proposition \ref{prop-key}, if WLP fails then those three generators must have a GCD of degree 1.  However, such a form is then automatically a socle element,  giving a contradiction.
\end{proof}

\begin{theorem} 
\[
e(r)=
\left \{
\begin{array}{ll}
+ \infty, & r \leq 2; \\
2 , & r = 3; \\
1, & r \geq 4.
\end{array}
\right.
\]

\end{theorem}

\begin{proof}
That $e(r)=+ \infty$ for $r\leq 2$, that is all level algebras of codimension at most 2 enjoy the WLP, is well-known (see \cite{HMNW, MZ1}). As for $r=3$, we know that there is a level (monomial) example without the WLP with Hilbert function $(1,3,5,5)$, constructed in \cite{zanello2}, Example 7. Thus, $e(3)\leq 2$. Hence, only level algebras with the following Hilbert functions need to be considered: $(1,3)$,   and $(1,3,a)$,  for $a=1,2,...,6$. 

A standard computation shows that applying Green's theorem and Lemma \ref{lll} takes care of all cases, except for $(1,3,3)$, for which we invoke Lemma \ref{133}.

Let $r\geq 4$. In light of Lemma \ref{222}, it remains to show that all level algebras $(1,r)$ enjoy the WLP, but this fact is trivial.
\end{proof}

\begin{proposition}
$e_2(3)=3.$
\end{proposition}

\begin{proof}
Notice that $e_2(3)\leq 3$, by Lemma \ref{brenner}. Thus it remains to show that all level algebras with the following Hilbert functions enjoy the WLP: $(1,3,2)$,   and $(1,3,a,2)$,  for $a=1,2,..., 6$. Note first of all that, for $a=1$ and 2, the set of such level algebras is empty (respectively, because of Macaulay's theorem and \cite{GHMS}, Proposition 3.8). The case $(1,3,3,2)$ follows from Lemma \ref{133} and the use of Green's theorem as in the previous proof (since $(2_{(3}))^{-1}_{-1}=0$). All the other cases are handled exactly as in the previous proof.
\end{proof}

\begin{remark} \label{hessian}
i) An entirely similar argument also proves that $e_3(3)=3$, that is the Brenner-Kaid example $(1,3,6,6,3)$ is the best possible in terms of socle degree for a codimension 3 level algebras of type 3 without the WLP.

ii) Let us now consider $e_1(r)$. In \cite{HMNW} Example 4.3, an example of a Gorenstein algebra with Hilbert function $(1,5,5,1)$ failing to have the WLP is provided. It is easy to extend that construction to non-WLP Gorenstein algebras with Hilbert function $(1,r,r,1)$ for all $r\geq 5$. This fact, combined with Green's theorem and Lemma \ref{lll}, easily implies that $e_1(r)=2$ for $r\geq 5$. 

As for the other values of $r$, we know that $e_1(r)=+\infty $ for $r=1,2$, and the results of the previous section show that $e_1(3)\geq 4$. We also asked whether $e_1(3)=+\infty $. As far as codimension 4 is concerned,  we have Ikeda's example (\cite{ikeda}, Example 4.4) of a non-WLP Gorenstein algebra with Hilbert function $(1,4,10,10,4,1)$. Furthermore, as pointed out to us by Junzo Watanabe, \cite{HMNW} Example 4.3 also shows that all Gorenstein algebras with Hilbert function $(1,4,4,1)$ enjoy the WLP (and actually more). Thus, $3\leq e_1(4)\leq 4$.

In fact, the referee of this paper pointed out to us that \cite{HMNW} Example 4.3 even implies the conclusion that $e_1(4) = 4$.  Indeed, Watanabe \cite{Wa} showed that  the Hessian of a form of degree $s$ is identically zero if and only if for the inverse system algebra $A$, the multiplication $\times L^{s-2} : A_1 \rightarrow A_{s-1}$ (where $L$ is a general linear form) does not have full rank.  Since Gordan and Noether \cite{GN} had shown that in four or fewer variables the vanishing of the Hessian implies that one of the variables can be eliminated, one can conclude that a Gorenstein algebra with $h$-vector $(1,4,a,4,1)$ has the Strong Lefschetz Property, and hence in particular WLP -- the bijectivity of the map from degree 1 to degree 3 implies the rest.  (A similar argument shows that for Gorenstein algebras with $h$-vectors $(1,3,n,n,\dots,n,n,3,1)$ or $(1,4,n,n,\dots,n,n,4,1)$, WLP implies SLP.)
\end{remark}

As promised, we now provide the examples of one level algebra of codimension 4 and type 2, and one of codimension 3 and type 4, without the WLP. We  omit the proofs, since they closely follow, respectively, that of \cite{zanello2}, Proposition 8, and that of Lemma \ref{brenner} of this paper (we used CoCoA \cite{CoCoA} for the computations).

\begin{example}\label{missing}
i) The codimension 4 and type 2 level algebra corresponding to the following inverse system module $M'\subset k[y_1,y_2,y_3,y_4]$, and having Hilbert function $(1,4,7,7,2)$, does not enjoy the WLP:
$$M'=\langle y_1^2y_2^2+y_1^2y_3^2+y_4^4,y_1^2y_2^2+y_2^2y_3^2+y_4^4\rangle .$$

ii) The codimension 3 and type 4 level algebra corresponding to the following inverse system module $M''\subset k[y_1,y_2,y_3]$, and having Hilbert function $(1,3,6,8,10,10,7,4)$, does not enjoy the WLP:
$$M''=\langle y_1^2y_3^5-y_1y_3^6,y_1^3y_3^4-y_1^5y_3^2,437y_1^7-232y_1^6y_2-423y_1^5y_2^2-567y_1^4y_2^3-$$$$-769y_1^3y_2^4+831y_1^2y_2^5-916y_1y_2^6-202y_2^7,(127y_1-548y_2-943y_3)^7\rangle .$$

\end{example} 
{\ }\\
{\bf Acknowledgements.} We warmly thank Uwe Nagel and Junzo Watanabe for helpful comments. We also thank the referee, for his or her insightful report.\\


\end{document}